\documentclass[12pt]{article}
\usepackage{amsmath,amsfonts,amssymb,amsthm,amscd}
\usepackage{color}
\usepackage{ulem}
\def\red{\textcolor{red}}

\title{Open subgroups of $p$-adic algebraic groups}
\date{\today}
\author{Anand Pillay\thanks{Supported by NSF grants DMS-1665035,   DMS-1760212 and DMS-2054271}\\University of Notre Dame \and Ningyuan Yao \thanks{Supported by National Social Science Fund of China, number 20CZX050 }\\Fudan University }
\newtheorem{Theorem}{Theorem}[section]
\newtheorem{Proposition}[Theorem]{Proposition}
\newtheorem{Definition}[Theorem]{Definition}
\newtheorem{Remark}[Theorem]{Remark}
\newtheorem{Lemma}[Theorem]{Lemma}
\newtheorem{Corollary}[Theorem]{Corollary}
\newtheorem{Fact}[Theorem]{Fact}

\newtheorem{Question}[Theorem]{Question}

\newtheorem*{Theorem-Add-I}{Additive Open Subgroup Theorem I}
\newtheorem*{Theorem-Add-II}{Additive Open Subgroup Theorem II}
\newtheorem*{Claim0}{Claim}

\newcommand{\R}{\mathbb R}   

\newcommand{\Q}{\mathbb Q}  
\newcommand{\Z}{\mathbb Z}  
\newcommand{\N}{\mathbb N}

\newcommand{\G}{\mathbb G}
\newcommand{\Qp}{{\mathbb Q}_p}
\newcommand{\Zp}{{\mathbb Z}_p}

\newcommand{\AI}{\operatorname{(Add-I)}}
\newcommand{\AII}{\operatorname{(Add-II)}}
\newcommand{\MI}{\operatorname{(Multi-I)}}
\newcommand{\MII}{\operatorname{(Multi-II)}}

\newcommand{\Ker}{\text{Ker}}
\begin{document}
\maketitle

\begin{abstract}
In \cite{Pillay} and more formally in \cite{Onshuus-Pillay} it was  asked whether open subgroups of $p$-adic algebraic groups are ($p$-adic) semialgebraic, equivalently, definable in the structure $(\Q_{p}, +, \times)$.  We give a positive answer in the commutative case. Together with results of \cite{Prasad} this leads to a positive answer for reductive algebraic groups.
\end{abstract}

\section{Introduction and preliminaries}
By a $p$-adic algebraic group we mean the group $G(\Q_{p})$ of $\Q_{p}$-points  of a (not necessarily linear) algebraic group $G$ defined over $\Q_{p}$. For our purposes we may assume $G$ to be connected as an algebraic group.  The group $G(\Q_{p})$ is a totally disconnected locally compact topological group with respect to the topology induced from the topological field
$\Q_{p}$.  The expressions {\em $p$-adic semialgebraic} and {\em first order definable (with parameters) in $(\Q_{p},+,\cdot)$} are synonymous and will be discussed later.

The following was asked, informally in \cite{Pillay}, and  more formally (Problem 1.1) in \cite{Onshuus-Pillay}.
\begin{Question} Let $G$ be an algebraic group over $\Q_{p}$. Is any open subgroup of $G(\Q_{p})$ $p$-adic semialgebraic?
\end{Question}

 In the $p$-adic topology, basic open neighbourhoods of a point $a\in \Q_{p}$ are of the form $\{x:v(x-a)\geq n\}$ where $n\in \Z$. These are $p$-adic semialgebraic, from which it follows that for any $p$-adic algebraic group $G(\Q_{p})$, there is a neighbourhood basis of the identity consisting of $p$-adic semialgebraic sets. We deduce from this  that any  {\em compact} open subgroup of any $G(\Q_{p})$ is semialgebraic.   When $G$ is semisimple and (almost) $\Q_{p}$-simple, then Question 1.1 has a positive answer.  This was pointed out to  the first author by Armand Borel in 1989 in a letter (basically a report on \cite{NP}).  It is based on a paper of Prasad  \cite{Prasad} which, among other things,  gives proofs of a (previously unpublished) theorem  of Rousseau and  Tits. The conclusion is that any open subgroup of $G(\Q_{p})$ is either semialgebraic and of finite index, or compact (so also semialgebraic).  Without being aware of Prasad's paper, in \cite{NP} an explicit proof was give in the case of $G = SL_{2}$, using the action of $SL_{2}(\Q_{p})$ on its tree, and we also gave a  positive answer to Question 1.1 for $GL_{2}(\Q_{p})$.

The main thrust of the current paper is to give a positive answer to Question 1.1 in the ``opposite" case when $G$ is a commutative (not necessarily linear) algebraic group over $\Q_{p}$.
We are also able to combine this with the earlier results to get for example a positive answer to Question 1.1 when $G$ is a reductive algebraic group over $\Q_{p}$.

Let us emphasize a fact which will be used a lot: if $G(\Q_{p})$ is a $p$-adic algebraic group, then $G$ has a neigborhood basis of the identity, consisting of compact open $p$-adic semialgebraic subgroups.

Question 1.1 was motivated   by the real case, where we DO have a positive answer: any open subgroup $H$ of a real algebraic group $G(\R)$  is (real) semialgebraic.  This is because, $G(\R)$ has finitely many semialgebraic, semialgebraically-connected components, each of which is actually connected. Hence the connected component $G(\R)^{0}$  of $G(\R)$ has finite index, is semialgebraic, and $G(\R)^{0}\leq H \leq G(\R)$, whereby $H$ is also semialgebraic.

Locally compact groups are built out of connected (real) Lie groups and totally disconnected groups, in the sense that if $G$ is locally compact, then its topological connected component $G^{0}$ is a projective limit of  connected Lie groups, and $G/G^{0}$ is   totally disconnected.  It is  somewhat interesting that in spite of the model-theoretic analogies between the fields of reals and $p$-adics (discussed below) the groups of the form $G(k)$, $G$ an algebraic group over $k$, for $k = \R$ or $\Q_{p}$ (more generally ``semialgebraic groups") live on the opposite sides of this behaviour: $G(\R)$ is a real Lie group, and $G(\Q_{p})$ is a $p$-adic Lie group, so totally disconnected.

Although model theory (in the form of definability) serves as some kind of motivation for  Question 1.1, it does not play any serious technical role in the proof.  Nevertheless we should give some background on definability and semialgebraicity. Our exposition and definitions are  tailor-made for the current paper.
The categories of real and $p$-adic semialgebraic sets can be given uniformly. Let $k$ be $\R$ or $\Q_{p}$.  Then the class of $k$-semialgebraic sets is the smallest class of subsets of Cartesian powers $k^{n}$ which (a) contains all $k$-algebraic sets, namely subsets of $k^{n}$ defined by finite systems of polynomial equations over $k$, and (b) is closed under (finite) Cartesian products and Boolean combinations and (c) if $m\leq n$, $\pi:k^{n}\to k^{m}$ is some coordinate projection, and $X\subseteq k^{n}$ is semialgebraic, then so is $\pi(X)$, and (d) If $\pi:k^{n}\to k^{m}$ is as in (c) and $X\subseteq  k^{m}$ is semialgebraic, so is $\pi^{-1}(X)$. So clearly the $k$-semialgebraic sets are  the sets definable with parameters in the structure $(k,+,\times)$.  This actually makes sense for any field $k$, although the expression ``semialgebraic" is typically reserved for the reals and $p$-adics.

It is well-known that the real semialgebraic subsets of $\R^{n}$ are precisely finite  unions of sets of the form  $\{{\bar a}\in \R^{n}: f_{1}({\bar a}) = 0 \wedge ... \wedge f_{r}({\bar a}) = 0 \wedge g({\bar a}) > 0\}$ where the $f_{i}$ and $g$ are polynomials over $\R$.  Note that  in $\R$, $x<y$ iff $x\neq y$ and $y-x$ is a square, and the right hand side {\em  is} semialgebraic as defined above.  Bearing this in mind, Macintyre's quantifier-elimination theorem for $\Q_{p}$ \cite{Macintyre} gives one of the strong model-theoretic  analogies between the reals and the $p$-adics, and says that any $p$-adic semialgebraic subset of $\Q_{p}^{n}$ is a finite Boolean combination of sets of the form  $\{{\bar a}: f({\bar a}) = 0\}$, and $\{{\bar a}: P_{m}(g({\bar a}))\}$  where $f,g$ are polynomials over $\Q_{p}$, $m=2,3,....$ and $P_{m}(x)$ stands for ``$x$ is an $mth$ power".
An important fact is that the valuation on $\Q_{p}$  is $p$-adic semialgebraic. Namely the ring of $p$-adic integers (the valuation ring) is definable (without parameters) by the formula
$\phi(x)$:  $\exists y(y^{2} = 1 + p^{3}x^{4})$. This works for any prime $p$.

We should also mention a related but distinct use of the expression ``semialgebraic" with respect to a valued field $(K,v)$. This  refers to subsets of $K^{n}$ which are Boolean combinations of $0$-sets of polynomials over $K$ and sets of the form $\{{\bar a}: v(f({\bar a})) \leq v(g({\bar a}))\}$ where $f,g$ are polynomials over $K$.  For example in the case where $(K,v)$ is a algebraically closed (nontrivially) valued field then the semialgebraic sets in this sense coincide with the sets definable in the structure $(K,+,\times, V)$ where $V$ is the valuation ring.

But in this paper we will be exclusively working with the valued field $\Q_{p}$ and semialgebraic has the earlier meaning.

We need information on the structure of connected commutative algebraic groups $G$ over a characteristic $0$ field $k$. It is convenient to refer to Milne's recent book \cite{Milne}.
First we have a unique connected linear algebraic subgroup $L$ of $G$ over $k$  such that $G/L$ is an abelian variety defined over $k$ (Theorem 8.27 of \cite{Milne}).  Secondly the linear part $L$ is a product of a vector group $V$ over $k$ and an algebraic torus $T$ over $k$.  (Corollary 16.15 of \cite{Milne}.) Where by a vector group over $k$ we mean an algebraic group over $k$ isomorphic over $k$ to a (finite) product of copies of the additive group ${\mathbb G}_{a}$. And by an algebraic torus over $k$ we mean an algebraic group over $k$ isomorphic over the {\em  algebraic closure}  of $k$ to a (finite) product of copies of the multiplicative group ${\mathbb G}_{m}$.  Given an algebraic torus $T$, we can write $T$ as an almost direct product of algebraic tori $T_{s}$ and $T_{a}$ over $k$, where $T_{s}$ (the split part of $T$)  is isomorphic over $k$ to a finite product of copies  of the multiplicative group ${\mathbb G}_{m}$, and $T_{a}$ is {\em anisotropic}, namely its character group over $k$ is trivial. (See Chapter 12 of \cite{Milne} and in particular Exercises 12-5 and 12-6.)

Now the underlying variety of an abelian variety $A$ over $k$ is a $k$-closed subvariety of some projective space ${\mathbb P}$. When $k = \Q_{p}$, ${\mathbb P}(k)$ is compact, as  is $A(k)$  (equipped witrh the $p$-adic topology).  And by Theorem 3.1 of \cite{Platonov-Rapinchuk}, again for $k = \Q_{p}$, and for $T$ an anisotropic torus over $k$, $T(k)$ is also compact.

For $k= \Qp$, let  $B=G/(V\times T_{s})$. Then we have an exact sequence  of algebraic groups over $k$, $1\to V\times T_{s} \to G \to B\to 1$, where the algebraic group $B$ is an extension of $A$ by $T_{a}$, and so $B(k)$ is also compact.   Putting all this together, we conclude:

\begin{Fact}
Let $G$ be a  commutative algebraic group over $\Q_{p}$ then there are one-dimensional connected linear algebraic groups $L_{1},..,L_{d}$ over $\Q_{p}$ and a connected algebraic group $B$ over $\Q_{p}$, and an exact sequence
$$1\to L = L_{1}\times ...\times L_{d} \to G \to B \to 1$$ over $\Q_{p}$.
Moreover,
\newline
(i) for each $i$, either $L_{i}$ is isomorphic over $\Q_{p}$ to  the additive group, or  isomorphic over $\Q_p$ to the multiplicative group,
\newline
(ii) $B(\Q_p)$ is compact, and the image of $G(\Q_{p})$ in $B(\Q_p)$ is an open (finite index) subgroup (so also compact).
\end{Fact}

The main result of this paper is:
\begin{Theorem}
 Let $G$ be a connected  commutative algebraic group over $\Q_{p}$. Then any open subgroup of $G(\Q_{p})$ is $p$-adic semialgebraic (definable).

\end{Theorem}

Let us remark that if $H$ is any abstract subgroup of  $G(\Q_{p})$ and $U$ is any open compact subgroup of $G(\Q_{p})$ then $H + U$ is an open subgroup of $G(\Q_{p})$.  This may suggest that applying this to the case where $G = \G_{m}\times \G_{a}$ and $H$ is the graph of $p$-adic exponentiation then we obtain non semialgebraic open subgroups of $G(\Q_{p})$.  But in fact the graph of  $p$-adic exponentiation lives in $\Z_{p}^{*} \times \Z_{p}$ so one does not really see it any more after adding a compact open subgroup.

\vspace{5mm}
\noindent
{\bf Acknowledgements.}  Thanks to Armand Borel and Guy Rousseau for their communications with the first author in 1989.
The material in Section 5 is influenced by what Borel told the first author in a letter dated 8 February 1989, and is also discussed in an Addendum to the paper \cite{NP}.
In a letter dated June 6th 1989 to the first author,  Rousseau  said again  how the result  in \cite{NP} that any proper open subgroup of $SL_{2}(\Q_{p})$ is  compact, is a special case of an unpublished result of Tits, with a proof also by Rousseau, and with a published account by Prasad.

The first author would also like to thank Dugald Macpherson for discussions in Autumn 2018 (during a meeting in Oaxaca).

\section{The proof: building blocks}
As already mentioned, if $G$ is an algebraic group over $\Q_{p}$ whose group of $\Q_p$-points is compact, then every open subgroup of $G(\Q_{p})$ is definable. We first use this to reduce the proof of the theorem to the case where $G$ is of the form $L_{1} \times ...\times L_{d}$ as in Fact 1.2.
\begin{Lemma}\label{Lemma-open-subgp-I}
Suppose $G$ is a commutative algebraic group over $\Q_{p}$, and $L\to G \to B$ is as in Fact 1.2. Suppose that $H$ is an open subgroup of $G(\Q_{p})$, and $H\cap L$ (an open subgroup of $L(\Q_p)$)  is definable, then $H$ is definable.
\end{Lemma}
\begin{proof}  Working with $\Q_p$ points we have, by Fact 1.2,  the induced short exact sequence:

$$   1\to L(Q_p) \to G_{1} \to_{\pi} B_1 \to 1  $$

where $G_{1} = G(\Q_p)$ and $B_{1}$ is an open finite index subgroup of  $B(\Q_p)$.

\vspace{2mm}
\noindent
Let $U$ be an open semialgebraic subgroup of $G_{1}$ which is contained in $H$.   Then $\pi(U)\leq \pi(H)$ are open subgroups of the compact group $B_{1}$, so $\pi(U)$ has finite index in  $\pi(H)$. Let $a_{1},..,a_{n}\in \pi(H)$ be representatives of the cosets of $\pi(U)$ in $\pi(H)$. Let $b_{i}\in H$ be such that $\pi(b_{i}) = a_{i}$.

\vspace{2mm}
\noindent
{\em Claim.}  $H = (H\cap L) + U +(b_{1}+U) .. + (b_{n} + U)$.
\newline
{\em Proof.}   Clearly the right hand side is contained in the left hand side.
Now suppose $h\in H$, so $\pi(h) = a_{i} + \pi(u)$ for some $i=1,..,n$ and some $u\in U$. But then $\pi(b_{i} + u) = \pi(h)$, so
As $b_{i} + u \in H$, $h-(b_{i} + u) \in  \Ker(\pi)  \cap H = H\cap L$, as required.                                                                   \qed

Hence as $H\cap L$ is definable, so is $H$.
\end{proof}

Hence, by Lemma 2.1  (and Fact 1.2) to prove Theorem 1.3, it suffices to prove:

\begin{Proposition}  Suppose $G =L_1\times...\times L_n$,    where each $L_{i}$ is either $\G_{a}$ of $\G_{m}$  (over $\Q_{p}$). Then  every open subgroup of $G(\Q_{p})$ is semialgebraic.
\end{Proposition}

Over the next few sections we will give the proof of Proposition 2.2, basically by induction on $d$. In  the remainder of this section we start with the case $d=1$ as well as some other observations.  $\Zp$ denotes the $p$-adic integers. We let $\Qp^{\times}$ denote the multiplicative group of the $p$-adic field, and $\Zp^{\times}$ its subgroup consisting of the invertible elements of $\Zp$ (i.e. elements of valuation $0$).

\begin{Lemma}  The proper open subgroups of $(\Qp,+)$ are precisely those defined by $v(x)\geq k$ for each $k\in \Z$, i.e. of the form $p^{k}\Zp$ for $k\in \Z$. In particular the proper open subgroups of $(\Qp,+)$ are all compact and semialgebraic.

\end{Lemma}
\begin{proof} We have mentioned in the introduction that these open balls $v(x)\geq k$ form a neighbourhood basis of the identity  in $(\Qp,+)$.
In particular if $H$ is an open subgroup of $(\Qp,+)$ then $H$ contains some $U =  ``v(x)\geq k"$ for some $k\in \Z$.  But the group $\Qp/U$ is isomorphic (abstractly) to the Pr\red{\"u}fer $p$-group $\Z(p^{\infty})$  (all complex $p^{n}th$ roots of unity for $n \geq 0$), and the only proper  abstract subgroups of $\Z(p^{\infty})$ are the subgroups of $p^{n}th$ roots of unity for each $n$. It follows that $H$ is of form $v(x)\geq k+n$ for the appropriate $n$.
\end{proof}

\begin{Lemma}\label{open subgp of Gm}  Let $H$ be an open subgroup  of $\Q_{p}^{\times}$. Then  either $H$ is compact, in which case  $(\Z_{p}^{\times})^{n} \leq H \leq \Z_{p}^{\times}$ for some $n$, or $H$ has a finite index subgroup of the form $(\Q_{p}^{\times})^{n}$ ($nth$ powers).  Either way $H$ is semialgebraic.
\end{Lemma}
\begin{proof} This is well-known,  and in fact appears explicitly in \cite{NP} (Fact 2.7).  Note that $H\cap \Z_{p}^{\times}$ is open of finite index. So if $H$ is contained in $\Z_{p}^{\times}$ then $(\Z_{p}^{\times})^{n} \leq H \leq \Z_{p}^{\times}$ for some $n$, and we know that $(\Z_{p}^{\times})^{n}$ has finite index in $\Z_{p}^{\times}$ so also in $H$.

Otherwise, $H\cdot \Z_{p}^{\times}/\Z_{p}^{\times}$ is a nontrivial subgroup of $\Q_{p}^{\times}/\Z_{p}^{\times} \cong  (\Z,+)$ (via the valuation $v$) so of the form $m\Z$ for some $m$, which has finite index in $\Z$.  As $H\cap \Z_{p}^{\times}$ has finite index in $\Z_{p}^{\times}$ it follows that $H$ has finite indexin $\Q_{p}^{\times}$ so contains some $(\Q_{p}^{\times})^{n}$ which has finite index in $\Q_{p}^{\times}$ so also has finite index in $H$.
\end{proof}



The following general and basic fact will be used a lot.

\begin{Lemma}\label{isomorphism}
Let $L_1 ,L_2$ be  groups and $G=L_1\times   L_2$. Suppose that $H$ is a  subgroup of  $G$.
Let $A_i=L_i\cap H$ be a normal subgroup of $L_i$ and $\pi_{i}$ the projection of $G_i$ on $L_i$, for $i=1,2$. Then there is an isomorphism $f_H:\pi_1(H)/A_1\rightarrow\pi_2(H)/A_2$ such that $H=\bigcup_{a\in \pi_1(H)}(aA_1)\times f_H(aA_1)$.
\end{Lemma}
\begin{proof}
For each $a\in \pi_1(H)$, there is $f_a\in \pi_2(H)$ such that $(a A_1)\times (f_a A_2)\subseteq H$. We claim that $f: a A_1\mapsto f_aA_2$ is  an isomorphism  from $\pi_1(G)/A_1$ to $\pi_2(G)/A_2$. First, if there are $b,c\in \pi_2(G)$ such that  $(a A_1)\times (b A_2)\subseteq H$ and $(a A_1)\times (c A_2)\subseteq H$, then $A_1\times (b^{-1}c A_2)\subseteq H$. We have  that $b^{-1}c A_2\subseteq H\cap L_2=A_2$, and thus $bA_2=cA_2$. So $f: a A_1\mapsto f_aA_2$ is well defined. Similarly, for each $b\in \pi_2(H)$, there is $g_b\in \pi_1(H)$ such that $(g_b A_1)\times (b A_2)\subseteq H$. We see that $g: b A_2\mapsto g_bA_1$ is the inverse of $f$. We conclude that  $f$ is a bijection. To see that $f$ is a morphism, take $a,b\in \pi_1(H)$, then \[
(aA_1\times f_aA_2)\cdot (bA_1\times f_bA_2)=(abA_1\times f_af_bA_2)=(abA_1\times f_{ab}A_2),
\]
which implies that $f$ is a morphism.

For each $(a,b)\in H$, we have $(a,b)(A_1\times A_2)=(aA_1\times bA_2)\subseteq H$. But we have showed that $bA_2=f_aA_2=f(aA_1)$, so $(a,b)\in aA_1\times f(aA_1)2$, hence $H\subseteq \bigcup_{a\in \pi_1(H)}(aA_1)\times f(aA_1)$ as required.
\end{proof}

\section{Product of Additive Groups}
Here we deal with Proposition 2.2 when all $L_{i}$ are the additive group.

\begin{Lemma}\label{Z(p-infty)-Automorphsiom}
Let $A$ be an open compact subgroup of $(\Qp,+)$ and  $\sigma: \Qp/A\rightarrow \Qp/A$ a group automorphism of $\Qp/A$. Then there is $\alpha\in \Zp^\times$ such that $\sigma(x+A)=\alpha x+A$ for all $x\in \Qp$.
\end{Lemma}
\begin{proof}
Let $O_n=(p^{-n}A)/A$. Then each $O_n$ is the cyclic group of order $p^n$, and the restriction of $\sigma$ to each $O_n$ is an automorphism of $O_n$. So there are integers $0<\alpha_1<\alpha_2<\alpha_3<...$ such that
\begin{itemize}
\item each $0<\alpha_n<p^{n}$ and
\item $\sigma(x)=\alpha_nx$ for $x\in O_n$.
\end{itemize}

\begin{Claim0}
$\alpha_n\equiv\alpha_{n+1}(\text{mod } p^{n})$ for each $n>0$.
\end{Claim0}
\begin{proof}
Induction on $n$. For $n=1$, we have that $O_1$ is the cyclic group of order $p$. So there is $0<\alpha_1<p$ such that $\sigma(  x)=\alpha_1   x$ for all $ x \in O_1$. Now suppose that  $\alpha_1<...<\alpha_n$ satisfy the claim. Let $x\in O_{n+1}$ be a generator of $O_{n+1}$, then $\sigma(x)=\alpha_{n+1}x$ and $\sigma(px)=p\sigma(x)=\alpha_{n+1}px$. As $px\in O_n$, we have $\sigma(px)=\alpha_npx$. So $\alpha_npx=\alpha_{n+1}px$. Since $x$ is a generator of $O_{n+1}$, we see that $p\alpha_n\equiv p\alpha_{n+1}(\text{mod } p^{n+1})$, so $\alpha_n\equiv\alpha_{n+1}(\text{mod } p^{n})$ as required.
\end{proof}
The above Claim shows that $\{\alpha_n|\ n\in\N^+\}$ is a Cauchy sequence in $\Zp^\times$, and  thus has a limit  $\alpha\in \Zp^\times$. Moreover, $\alpha\equiv\alpha_{n}(\text{mod } p^{n})$ for each $n\in\N^+$, which means that $\alpha x+A=\alpha_{n} x+A$ for all $x\in p^{-n}A$. We conclude that for each $x\in p^{-n}A$, $\sigma(x+A)=\alpha_n(x+A)=\alpha_nx+A=\alpha x+A$. This completes the proof.
\end{proof}

\begin{Corollary}\label{Qp/A-to-Qp/B}
Let $A$ and $B$ be  open compact subgroups of $(\Qp,+)$ and  $\eta: \Qp/A\rightarrow \Qp/B$ a group isomorphism. Then there is $\alpha\in \Qp^{\times}$ such that $\eta(x+A)=\beta x+A$ for all $x\in \Qp$.
\end{Corollary}
\begin{proof}
Suppose that $A=p^n\Zp$ and $B=p^m\Zp$. Let $f:\Qp\rightarrow \Qp$ be the map defined by $ x\mapsto p^{m-n}x$. Then $\bar f: x/A\mapsto f(x)/B$ is an  isomorphism from $\Qp/A$ to $\Qp/B$. As ${\bar{f}}^{-1}\circ\eta$ is an automorphism of $\Qp/A$,  we see from Lemma \ref{Z(p-infty)-Automorphsiom} that $ {\bar{f}}^{-1}\circ\eta(x)=\alpha x+A$ for some $\alpha\in\Zp^\times$, So $\eta(x)=\bar f\circ ({\bar{f}}^{-1}\circ\eta)(x)=p^{m-n}\alpha+A$.
\end{proof}


\begin{Theorem}[Additive Open Subgroup Theorem I]\label{Additive-Open-Subgroup-is-definable}
Let $G=L_1\times...\times L_n$  where each $L_{i} = (\Q_{p},+)$. Let  $H$ an open subgroup of $G$ such that $H_A=H\cap L_1$ is compact and $\pi_1(H)=L_1$. Let $H_B=H\cap (L_2\times ..\times L_n)$. There there is  a definable linear map $f:L_1 \rightarrow L_2\times...\times L_n$ such that
\[
H=\bigcup_{a\in \Q_{p}}(a+H_A)\times(f(a)+H_B ).
\]
\end{Theorem}

\begin{Theorem}[Additive Open Subgroup Theorem II]
Suppose that $H$ is an   open subgroup of $G=L_1\times...\times L_n$ where again each $L_{i}$ is $(\Q_{p},+)$. If $O$ is an   open subgroup of $H$ such that $H/O\cong \Z(p^\infty)$, then   there is a definable linear map $h: \Qp\rightarrow H$  inducing an isomorphism between $\Q_{p}/\Z_{p}$ and  $H/O$.
\end{Theorem}

$\AI_{n}$ will denote the statement of the Additive Open Subgroup Theorem for the case $n$. Likewise for $\AII_{n}$.


\begin{Lemma}\label{ IIn-imp-I(n+1)}
$\AI_n, \AII_n\implies \AI_{n+1}$ for $n\geq 1$
\end{Lemma}

\begin{proof}
 Assume that  $G=L_1\times...\times L_{n+1}$. Let $H$, $H_A$ and $H_B$ be as in the {Additive Open Subgroup Theorem I}. Let $\pi_{(2,...,n+1)}$ be the projection from $G$ onto $L_2\times...\times L_{n+1}$.   By Lemma \ref{isomorphism}, there is a map $g:L_1\rightarrow \pi_{(2,...,n+1)}(H)$ such that
 \[
 \bar g: L_1/H_A\rightarrow \pi_{(2,...,n+1)}(H)/H_B, \ a/H_A\mapsto g(a)/H_B
 \]
 is an isomorphism and
 \[
 H=\bigcup_{a\in\Qp}(a+H_A)\times (g(a)+H_B).
 \]
 Since $\pi_{(2,...,n+1)}(H)/H_B\cong L_1/H_A\cong \Z(p^\infty)$, applying $\AII_n$, there is a definable linear map $h:\Qp\rightarrow \pi_{(2,...,n+1)}(H_B)$ such that $\bar h: L_1/H_A\rightarrow \pi_{(2,...,n+1)}(H)/H_B$ is an isomorphism. Now $\bar h$ and $\bar g$ induce an automorphism $\bar \sigma$ of $L_1/H_A$ such that $  \bar g =\bar h\circ \bar \sigma$. By Corollary \ref{Qp/A-to-Qp/B}, there is $\gamma\in\Zp^\times$ such that $\bar \sigma(a/H_A)=\gamma a/H_A$. Let $\sigma: \Qp\rightarrow\Qp$ defined by $x\mapsto \gamma x$ and $f=h\circ \sigma$, then $f$ is definable linear and $f(x)/H_B= g(x)/H_B$ for all $x\in\Qp$. This completes the proof.
\end{proof}

\begin{Lemma}\label{I(n+1)-and-IIn-imp-II(n+1)}
$\AI_{n+1},\AII_n \implies \AII_{n+1} $ where $n\geq 1$.
\end{Lemma}
\begin{proof}
Assume that  $G=L_1\times...\times L_{n+1}$. Let $H$, $H_A$, and $H_B$ be as in the {Additive Open Subgroup Theorem I}.
Let $\pi_{(2,...n+1)}$ be the projection from $G$ onto $L_2\times...\times L_{n+1}$,  $O_A=O\cap L_1$, and  $O_B=O\cap L_2\times...\times L_n$.

\vspace{0.5cm}
\noindent $\bullet$  {Case 1}:  $H_A=O_A=L_1$. Then $H=L_1\times H_B$ and $O=L_1\times O_B$. We see that  $H/O\cong H_B/O_B$   and we can apply $\AII_{n}$ to $H_B/O_B$.

\vspace{0.5cm}
\noindent  $\bullet$  {Case 2}:  $H_A=L_1$, $\pi_1(O)=L_1$, and $O_A$ is compact.   Then $H=H_A\times H_B$. 

Let ${h_0}: \Qp\rightarrow H$ defined by $a\mapsto (a,0,...,0)$. Then  it is easy to see that $\bar {h_0}: a/O_A\mapsto h_0(a)/O$ is a group homomorphism from $\Qp/O_A$ to $H/O$.
On the other side,  $a/O_A\in \Qp/O$ has order $p^m$ iff $  p^{m}a\in O_A$ and $ p^{-(m-1)}a\notin O_A$, which means that  $p^m (a,0,...,0) \in O_A\times O_B\subseteq O$ and $p^{m-1} (a,0,...,0)\in (p^{m-1}a+O_A)\times O_B$. Since $((p^{m-1}a+O_A)\times O_B)\cap O=\emptyset$, we have that $(a,0,...,0)$ has order $p^m$ in $H/O$. We see that  $\bar h_0$  is injective and the image of   $\bar h_0$  is infinite. So $\bar {h_0}$ is an isomorphism. Let $h_1: \Qp\rightarrow H$ be a definable linear function such that $a/\Zp\mapsto h_1(a)/O_A$ is an isomorphism from $\Qp/\Zp$ to $\Qp/O_A$. Let $h=h_1\circ h_0$, then $h$ is definable linear,  thus satisfies the property.

\vspace{0.5cm}
\noindent $\bullet$  {Case 3}:   $\pi_1(H)$ is compact. Then $H=\bigcup_{a\in \pi_1(H)}(a+H_A)\times (g(a)+H_B)$ and $H/(H_A\times H_B)$ is finite.
Consider the short exact sequence
\[
1 \rightarrow (H_B+O)/O\rightarrow H/O\rightarrow H/(H_B+O)\rightarrow 1
\]
As $O_A\times H_B\le H_B+O$,  we see that $(H_A\times H_B)\cap (H_B+O)$ is a finite-index subgroup of $H_A\times H_B$, so $H/(H_B+O)$ is finite. Since $H/O = \Z(p^\infty)$ has no proper infinite subgroup, we have that $(H_B+O)=H$. Now
\[
(H_B+O)/O\cong H_B/(H_B\cap O) \text{ via } (0,a)/O\mapsto a/(H_B\cap O).
\]
Applying $\AII_n$ to $H_B/(H_B\cap O)$, there is a definable linear function $h_0:\Qp\rightarrow H_B$ such that $a/\Zp\mapsto h_0(a)/H_B\cap O$ is an isomorphism. So
\[
h:\Qp\rightarrow H,\ \  x\mapsto (0,h_0(x))
\]
satisfies the property.

\vspace{0.5cm}
\noindent
$\bullet$ {Case 4}:   $H_A$ and $\pi_1(O)$ are compact and $\pi_1(H)=L_1$. By $\AI_{n+1}$, there is a definable linear function $f:\Qp\rightarrow L_2\times ...\times L_{n+1}$ such that
\[
H=\bigcup_{a\in \Qp}(a+H_A)\times (f(a)+H_B).
\]

Let $D=H\cap (\pi_1(O)\times L_2\times...\times L_n)$. Then $O\subseteq D$. Since $H/D$ is infinite, we see that $D/O$ is a finite subgroup of $H/O$. Suppose that $D/O$ has order $p^n$, then for each $d\in D$, $p^nd\in O$ and thus
\[
D=p^{-n}O=\{p^{-n}x|\ x\in O\}.
\]
Consider the map $h: \Qp\rightarrow H$, $a\mapsto ( a, f(a))$. We claim that
\[
\bar h: \Qp/p^{n}\pi_1(O)\rightarrow H/O,\ \  a/p^{n}\pi_1(O)\mapsto h(a)/O
\]
is a group morphism. To see   $\bar h$ is a group morphism, we only need to check that $h(a)\in O$ for each $a\in p^{n}\pi_1(O)$. Take $a\in p^n\pi_1(O)$, then $a=p^nb$ with $b\in \pi_1(O)$. So $(b,f(b))\in D$ and  thus $p^n(b,f(b))= (p^n b,p^nf(b)) \in O$. Since $f$ is linear, we have $p^nf(b)=f(p^nb)=f(a)$. So $h(a)=(a,f(a))\in O$ as required.

Let $a\in \Qp$ such that $a/ p^n\pi_1(O)$ has order $p^m$ in $\Qp/p^n\pi_1(O)$, then $a=p^{n-m}b$  for some  $b\in \pi_1(O)\backslash p\pi_1(O)$.  So $(b,f(b))\in D\backslash pD$. which implies that $(b,f(b))/O$ has order $p^n$ in $D/O$, and thus
\[
p^{n-m}(b,f(b))/O=(a,f(a))/O=h(a)/O
\]
has order $p^m$. We conclude that $\bar h$ is injective, so is isomorphsim.

\vspace{0.4cm}
\noindent
$\bullet$ {Case 5}:  $O_A\leq H_A$ is compact and $\pi_1(H)=\pi_1(O)=L_1$. Assume that
\[
H=\bigcup_{a\in \Qp}(a+H_A)\times (f_H(a)+H_B)
\]
and
\[
O=\bigcup_{a\in \Qp}(a+O_A)\times (f_O(a)+O_B)
\]
where $f_H$ and $f_O$ are functions from $\Qp$ to $L_2\times ...\times L_{n+1}$.
It is easy to see that $f_O(a)-f_H(a)\in H_B$ for each $a\in\Qp$. So $H=(H_A\times H_B)+O)$, and thus $(H_A\times H_B)/((H_A\times H_B)\cap O)$ is isomorphic to $H/O=(H_A\times H_B)+O)/O$   via
\[
x/((H_A\times H_B)\cap O)\mapsto x/O
\]
for $x\in H_A\times H_B$. By Case 3, we see that there is a definable linear function $h:L_1\rightarrow H_A\times H_B$ such that $a/\Zp\mapsto h(a)/((H_A\times H_B)\cap O)$ is an isomorphism from $\Qp/\Zp$ to $(H_A\times H_B)/((H_A\times H_B)\cap O)$, so $a/\Zp\mapsto h(a)/O$ is also an isomorphism from $\Qp/\Zp$ to $H/O$.
\end{proof}

It is easy to see from Corollary \ref{Qp/A-to-Qp/B} that $\AII_1$ holds. Using Lemma \ref{ IIn-imp-I(n+1)} and Lemma \ref{I(n+1)-and-IIn-imp-II(n+1)}, we can proof Theorem \ref{Additive-Open-Subgroup-is-definable} by induction on $n\in\N^+$.

\begin{Theorem}
Let $H$ be an open subgroup of $G=L_1\times...\times L_n$, where each $L_i= (\Q_{p},+)$. Then $H$ is definable.
\end{Theorem}
\begin{proof}
Induction on $n$. Let $H\cap L_1=H_A$ and $H_B=H\cap (L_2\times...\times L_n)$. Then $H_A$ and $H_B$ are definable by induction hypothesis.
If $\pi_1(H)/H_A$ is finite, then $H_A\times H_B$ is a finite index subgroup of $H$, thus $H$ is definable.
If $\pi_1(H)/H_A$ is infinite, then $\pi_1(H)=L_1$ and $H_A$ is compact. By Theorem \ref{Additive-Open-Subgroup-is-definable}, there is a definable linear function $f: L_1\rightarrow L_2\times...\times L_n$ such that
\[
H=\bigcup_{a\in \Qp}(a+H_A)\times (f(a)+H_B).
\]
For any $a\in L_1$ and $b\in L_2\times...\times L_n$, $(a,b)\in H\iff b-f(a)\in H_B$. So $H$ is definable.
\end{proof}

\begin{Lemma}
Suppose that $G=L_1\times...\times L_n$, with $L_i= (\Q_{p},+)$ for each $i\leq n$. Let  $H$ be an open subgroup of $G(\Qp)$, then $G/H$ is a torsion group.
\end{Lemma}
\begin{proof}
Induction on $n$. For $n=1$, $G/H$ is either trivial or isomorphic to  $\Z (p^\infty)$. Now assume that $G=L_1\times...\times L_{n+1}$. Let $G_1=L_1\times...\times L_n$, then we have a short exact sequence
\[
1\rightarrow G_1/G_1\cap H\rightarrow G/H\rightarrow L_{n+1}/\pi_{n+1} (H)\rightarrow 1,
\]
 where $\pi_{n+1}$ is the projection from $G$ onto $L_{n+1}$.  We see from induction hypothesis that both $G_1/G_1\cap H$ and $L_{n+1}/ \pi_{n+1}  (H)$ are torsion groups. So $G/H$ is also a torsion group.
\end{proof}

\begin{Lemma}\label{Quotient-group-of-add-group}
Suppose that $G=L_1\times...\times L_n$, with $L_i= (\Q_{p},+)$ for each $i\leq n$. Let  $O\leq H$ be   open subgroups of $G$, then $H/O$ is a torsion group.
\end{Lemma}
\begin{proof}
 $H/O$ is a subgroup of $G/O$, so is a torsion group.
\end{proof}

\section{Product of Multiplicative Groups}

As above $\Q_{p}^{\times}$ denotes the multiplicative group of $\Q_{p}$. Likewise for $\Z_{p}^{\times}$.   $(\Q_{p}^{\times})^{k}$ denotes the $kth$ powers of the multiplicative group.

\begin{Theorem}[{Multiplicative Open Subgroup Theorem I}]\label{Multi-Open Subgroup Theorem I}
Let $G=L_1\times...\times L_n$ with each $L_i= \Q_{p}^{\times}$, and  $H$ an open subgroup of $G$. Then $H$ is definable.
\end{Theorem}

\begin{Theorem}[{Multiplicative Open Subgroup Theorem II}]\label{Multi-Open-Subgroup II}
Suppose that $H$ is an   open subgroup of $G=L_1\times...\times L_n$ with each $L_i=\Q_{p}^{\times}$. If $O$ is an   open subgroup of $H$ such that $H/O\cong (\Z,+)$, then   there is a definable map   $h: \Qp\rightarrow H$ such that $a/\Zp^\times\mapsto h(a)/O$ is an isomorphism from $\Q_{p}^{\times}/\Zp^\times$ to $H/O$.
\end{Theorem}
The  Multiplicative Open Subgroup Theorem I (resp. II) for  $n$ will be denoted by $\MI_n$ (resp. $\MII_n$).  Note that $\MI_1$ holds by Lemma \ref{open subgp of Gm}.

\begin{Remark}
Since $\text{Th}(\Qp)$ has definable Skolem functions (see \cite{definable Skolem functions}), we have the following:  If $H$ is a definable subgroup of $L_1\times...\times L_{n}$ (where $L_{i}$ is $\Q_{p}^{\times}$)
and $H_A=H\cap L_1$ and $H_B=H\cap (L_2\times ..\times L_n)$, then there is the a definable map
\[
f:\pi_1(H)\rightarrow L_2\times...\times L_n
\]
such that
\[
H=\bigcup_{a\in \pi_1(H)}(a H_A)\times(f(a) H_B ).
\]
\end{Remark}

\begin{Remark}\label{rmk-Pk}
 Note that for each positive integer $k$,
\[
f_k: \Qp^\times \to (\Q_{p}^{\times})^{k}, \ x\mapsto x^k
\]
is a definable map such that $\bar f_{k}: x/\Zp^\times\mapsto f_k(x)/\Zp^\times$ is an isomorphism from $\Q_{p}^{\times}/\Zp^\times$ to $(\Q_{p}^{\times})^{k}\Zp^{\times}/\Zp^\times$. Since $\text{Th}(\Qp)$ has definable Skolem functions, $f_k$ has a definable section $h_k:   (\Q_{p}^{\times})^{k} \rightarrow \Q_{p}^{\times}$. Clearly, $\bar {h_k}: x/\Zp^\times\mapsto h_k(x)/\Zp^\times$ is an isomorphism from $(\Q_{p}^{\times})^{k}\Zp^{\times}/\Zp^\times$ to $\Q_{p}^{\times}/\Zp^\times$.
\end{Remark}

\begin{Lemma}\label{MII1}
 $ \MII_1$ holds.
\end{Lemma}
\begin{proof}
Let $O\leq H$ be open subgroups of $\Q_{p}^{\times}$ such that $H/O\cong(\Z,+)$. Then $O$ is compact and thus a subgroup of $\Zp^\times$. Let $a\in H$ have minimal positive valuation. Then $H/O=\{a^m/O|\ m\in \Z\}$. Suppose that $v(a)=d$, and let
\[
f:\Q_{p}^{\times}\rightarrow \Q_{p}^{\times},\ \ x\mapsto x^d,
\]
then
\[
\bar f: \Q_{p}^{\times}/\Zp^\times\rightarrow H/O,\ \  x/\Zp^\times \mapsto x^{d}/O
\]
is an isomorphism.
\end{proof}

\begin{Lemma}\label{MIIn+1}
 $\MI_{n+1},\ \MII_n\implies \MII_{n+1}$  for $n\geq 1$.
\end{Lemma}

\begin{proof}
Assume that  $G=L_1\times...\times L_{n+1}$ where each $L_i=\Q_{p}^{\times}$. Let $O\leq H$ be  open subgroups of $G$, $H_A=H\cap L_1$,   $H_B=H\cap (L_2\times ..\times L_{n+1})$,  $O_A=O\cap L_1$, and  $O_B=O\cap (L_2\times...\times L_{n+1})$.

\vspace{0.4cm}
\noindent $\bullet$ {Case 1}:  $H_A$ and $O_A$ have finite index in $\pi_1(H)$. Then $H_A\times H_B$ has finite index in $H$ and  $O_A\times O_B$ has finite index in $O$.

Consider the short exact sequence
\[
1\rightarrow (H_A\times H_B)O/O\rightarrow H/O\rightarrow H/(H_A\times H_B)O\rightarrow 1.
\]

Since $H/H_A\times H_B$ is finite, we see that $H/((H_A\times H_B)O)$ is finite, so $(H_A\times H_B)O/O$ is an infinite subgroup of $H/O$, and thus is isomorphic to $(\Z,+)$. Assume that $|H/((H_A\times H_B)O)|=k$. Note that the subgroup of $(\Z,+)$ of index $k$ is of the form $(k\Z,+)$.
Let $b\in H$ such that $b/O$ is the generator of $H/O$, then
\[
H=\bigcup_{n\in \mathbb Z}b^nO  \ \ \text{and} \ \ (H_A\times H_B)O=\bigcup_{n\in \mathbb Z}b^{nk}O,
\]
which means that
\[
H=\bigcup_{i=0}^{k-1}b^i(H_A\times H_B)O.
\]

Since
\[
H_AO/O\cong H_A/(H_A\cap O)=H_A/O_A
\]
is a finite subgroup of $H/O$, we see that $H_A=O_A$ as $(\Z,+)$ has no non-trivial finite subgroups.  So
\[
(H_A\times H_B)O/O\cong (H_A\times H_B)/((H_A\times H_B)\cap O)\cong H_B/O_B.
\]
Applying $\MII_n$ to $H_B/O_B$, and by Remark \ref{rmk-Pk}, there is a definable map   $h_0: (\Q_{p}^{\times})^{k}\rightarrow (H_A\times H_B)O$  such that
\[
\bar h_0: (\Q_{p}^{\times})^{k}\Zp^\times /\Zp^\times\rightarrow (H_A\times H_B)O/O,\ a/\Zp^\times \mapsto h_0(a)/O
\]
is an isomorphism. Since $p^k/\Zp^\times$ is the generator of $(\Q_{p}^{\times})^{k}\Zp^\times /\Zp^\times$, we have $h_0(p^k)/O$ is the generator of $(H_A\times H_B)O/O$. So $h_0(p^k)/O=b^k/O$. For any $x\in (\Q_{p}^{\times})^{k}$, if $v(x)= mk$, then $x/\Zp^\times = {(p^k)}^m/\Zp^\times$, so
\[
h_0(x)/O=h_0({(p^k)}^m)/O={(b^k)}^m/O=b^{km}/O=b^{v(x)}/O.
\]

Now $\Q_{p}^{\times}=\bigcup_{i=0}^{k-1}p^i\Zp^\times (\Q_{p}^{\times})^{k}$. Applying definable Sklolem functions, there is a definable function $d:\Q_{p}^{\times}\rightarrow (\Q_{p}^{\times})^{k}$  such that $xd(x)^{-1}=p^i\Zp^\times$ if $x\in p^i\Zp^\times (\Q_{p}^{\times})^{k}$ for $i=0,...,k-1$. Let $h: \Q_{p}^{\times}\rightarrow H$ defined by $h(x)=b^ih_0(d(x))$ if $x\in p^i\Zp^\times (\Q_{p}^{\times})^{k}$ for  $i=0,...,k-1$. Then $h$ is a definable function.

To see that
\[
\bar h:\Q_{p}^{\times}/\Zp^\times\rightarrow H/O,\ \  x/\Z_p^\times\mapsto  h(x)/O
\]
is an isomorphism , we only need to check that $h(x)/O=b^{v(x)}/O$ for all $x\in \Q_{p}^{\times}$. Let $x\in\Q_{p}^{\times}$. Take $0\leq i<k$ such that $v(x)=v(d(x))+i$, then \[
h(x)/O=b^ih_0(d(x))/O=b^ib^{v(d(x))}/O=b^{v(x)}/O.
\]
This completes the proof the Case 1.

\vspace{0.4cm}
\noindent  $\bullet$ {Case 2}:  $H_A$ is non-compact and  $O_A$ is compact.  Then $H_AO/O\cong H_A/O_A$ is an infinite subgroup of $H/O$ and thus $H/(H_AO)$ is finite (using our assumption that $H/O$ is isomorphic to $(\Z,+)$). Assume that $|H/(H_AO)|=k$.

Applying $\MII_n$ to $H_A/O_A$, there is a definable map   $h_0: (\Q_{p}^{\times})^{k}\rightarrow H_AO$  such that
\[
\bar h_0: (\Q_{p}^{\times})^{k}\Zp^\times /\Zp^\times\rightarrow H_AO/O,\ a/\Zp^\times \mapsto h_0(a)/O
\]
is an isomorphism. Let $d:\Q_{p}^{\times}\rightarrow (\Q_{p}^{\times})^{k}$ be a definable function that $xd(x)^{-1}=p^i\Zp^\times$ if $x\in p^i\Zp^\times (\Q_{p}^{\times})^{k}$ for $i=0,...,k-1$. Let $h: \Q_{p}^{\times}\rightarrow H$ defined by $h(x)=b^ih_0(d(x))$ if $x\in p^i\Zp^\times (\Q^{\times})^{k}$ for  $i=0,...,k-1$. Then a similar argument as in Case 1 shows that $h$ is a definable function and
\[
\bar h:\Q_{p}^{\times}/\Zp^\times\rightarrow H/O,\ \  x/\Z_p^\times\mapsto  h(x)/O
\]
is an isomorphism.

\vspace{0.4cm}
\noindent
$\bullet$ {Case 3}:    $H_A$ and $\pi_1(O)$ are compact and $\pi_1(H)$ is non-compact.
Consider the short exact sequence
\[
1\rightarrow (H_A\times H_B)O/O\rightarrow H/O\rightarrow H/(H_A\times H_B)O\rightarrow 1.
\]
Since $\pi_1(H)$ is non-compact and $\pi((H_A\times H_B)O)= H_A\pi_1(O)$ is compact, we see that $H/(H_A\times H_B)O$ is infinite. Since the quotient groups of $(\Z,+)$ are either finite or $(\Z,+)$ itself, we have  $H_A\times H_B\leq O$. Consider the short exact sequences
\[
1\to O/(H_A\times H_B)\to H/(H_A\times H_B)\to H/O\to 1.
\]
and
\[
1\to \pi_1(O)/H_A \to \pi_1(H)/H_A\to \pi_1(H)/\pi_1(O)\to 1.
\]
Since $O/(H_A\times H_B)\cong \pi_1(O)/H_A$ via $x/(H_A\times H_B)\mapsto \pi_1(x)/H_A$, and $H/(H_A\times H_B)\cong \pi_1(H)/H_A$ also via $x/(H_A\times H_B)\mapsto \pi_1(x)/H_A$, we see that
\[
\bar \pi_1: H/O\to  \pi_1(H)/\pi_1(O),\  x/O\mapsto \pi_1(x)/\pi_1(O)
\]
is an isomorphism.  Let $h_0$ be a definable section of $\pi_1$, then
\[
\bar h_0:  \pi_1(H)/\pi_1(O)\rightarrow H/O,\ a/\pi(O) \mapsto h_0(a)/O.
\]
is an isomorphism. Applying $\MII_1$ to $\pi_1(H)/\pi_1(O)$, there is a definable map   $h_1: \Q_{p}^{\times}\rightarrow  \pi_1(H)$  such that
\[
\bar h_1:  \Q_{p}^{\times}   /\Zp^\times\rightarrow  \pi_1(H)/\pi_1(O),\ a/\Zp^\times \mapsto h_1(a)/\pi_1(O).
\]
is an isomorphism.
Applying $\MI_{n+1}$ to $H$, there is definable map
\[
f:\pi_1(H)\rightarrow L_2\times...\times L_{n+1}
\]
such that
\[
H=\bigcup_{a\in\pi_1(H)}(aH_A)\times(f(a)H_B).
\]
Let
\[
h: \Q_{p}^{\times}\rightarrow H, x\mapsto (h_0(h_1(x)),f(h_0(h_1(x)))),
\]
Then
\[
\bar h :  \Q_{p}^{\times}  /\Zp^\times\rightarrow H/O,\ a/\Zp^\times \mapsto h(a)/H_A.
\]
is an isomorphism. This completes the proof of Case 3.

\vspace{0.4cm}
\noindent
$\bullet$ {Case 4}:  $ H_A$ is compact and $\pi_1(O)$ is non-compact. Applying $\MI_{n+1}$ to $H$ and $O$, let $f_H$ and $f_O$ be definable functions from $L_1$ to $L_2\times ...\times L_{n+1}$ such that
\[
H=\bigcup_{a\in \pi_1(H)}(aH_A)\times (f_H(a)H_B)
\]
and
\[
O=\bigcup_{a\in \pi_1(O)}(aO_A)\times (f_O(a)O_B)
\]

It is easy to see that $f_O(a)-f_H(a)\in H_B$ for each $a\in\pi_1(O)$. So
\[
 (H_A\times H_B)O=H\cap (\pi_1(O)H_A\times L_2\times...\times L_{n+1}),
\]
and thus  $H/(H_A\times H_B)O\cong \pi_1(H)/\pi_1(O)H_A$ is finite. We see that $(H_A\times H_B)O/O$ is an infinite subgroup of $H/O$ and hence is isomorphic to $\Z$. Since $H_A$ is compact,
\[
H_AO/O\cong H_A/O\cap H_A=H_A/O_A
\]
is a finite subgroup of $(H_A\times H_B)O/O$, thus is trivial, we conclude that $H_A=O_A$, so
\[
(H_A\times H_B)O/O\cong (H_A\times H_B)/(O\cap H_A\times H_B)=(H_A\times H_B)/(H_A\times O_B)\cong H_B/O_B.
\]
Assume that $|H/(H_A\times H_B)O|=k$. Applying $\MII_n$ to $H_B/O_B$, there is a definable map   $h_0: (\Q_{p}^{\times})^{k}\rightarrow (H_A\times H_B)O $  such that
\[
\bar h_0 :(\Q_{p}^{\times})^{k} \Zp^\times /\Zp^\times\rightarrow (H_A\times H_B)O/O,\ \ a/\Zp^\times \mapsto h_0(a)/O
\]
is an isomorphism.

Let $d:\Q_{p}^{\times}\rightarrow (\Q_{p}^{\times})^{k}$ be a definable function  such  that $xd(x)^{-1}=p^i\Zp^\times$ if $x\in p^i\Zp^\times(\Q_{p}^{\times})^{k}$ for $i=0,...,k-1$. Let $h: \Q_{p}^{\times}\rightarrow H$ be defined by $h(x)=b^ih_0(d(x))$ if $x\in p^i\Zp^\times (\Q_{p}^{\times})^{k}$ for  $i=0,...,k-1$. Then a similar argument as in Case 1 shows that $h$ is a definable function and
\[
\bar h:\Q_{p}^{\times}/\Zp^\times\rightarrow H/O,\ \ x/\Z_p^\times\mapsto  h(x)/O
\]
is an isomorphism.

\end{proof}

\begin{Lemma}\label{MIn+1}
$\MI_n,\ \MII_n\implies \MI_{n+1}$ for  $n\geq 1$.
\end{Lemma}
\begin{proof}
Let $H$ be an open subgroup of $L_1\times...\times L_{n+1}$. Let $H_A$ and $H_B$ be as the above. If $\pi_1(H)/H_A$ is finite, then $H/H_A\times H_B$ is finite. By $\MI_n$, $H_A$ and $H_B$ are definable, so is $H$.

If $\pi_1(H)/H_A$ is infinite, then $H_A\subseteq \Zp^\times$, and $\pi_1(H)$ is non-compact. Take a function  $f:\pi_1(H)\rightarrow L_2\times...\times L_{n+1}$  such that
\[
H=\bigcup_{a\in \pi_1(H)}(aH_A)\times (f(a)H_B).
\]
Let $C_A =\pi_1(H)\cap \Zp^\times$. Consider the short exact sequence
\[
1\rightarrow C_A/H_A\rightarrow \pi_1(H)/H_A \rightarrow \pi_1(H)/C_A\rightarrow 1.
\]
As $\pi_1(H)/C_A$ is isomorphic to $\Z$, we see that $C_A/H_A$ is precisely the subgroup of torsion elements of $ \pi_1(H)/H_A $.
Take  $\eta\in \pi_1(H)$ with minimal positive valuation. Let $\theta_{1},...,\theta_{k}\in C_A$ such that $C_A$ is a disjoint union of $\{\theta_{1} H_A,...,\theta_{k} H_A\}$. Then $\pi_1(H)$ is a disjoint union of $\{\theta_{i}\eta^{m}H_A|\ 1\leq i\leq k,\ m\in\Z\}$. Let $f(\theta_i)=\gamma_i$ for $1\leq i\leq k$, then
 \[
 H=\bigcup_{1\leq i\leq k,m\in\Z}(\theta_{i}\eta^{m}H_A)\times  (\gamma_{i}f(\eta^{m})H_B).
 \]
Let $H^*=HC_A$, then $\pi_1(H^*) =\pi_1(H)$ and $\pi_{(2,...n+1)}(H^*)=\pi_{(2,...n+1)}(H)$.
As $\pi_1(H^*)\cap \Zp^\times =C_A$, we have   that $\pi_1(H^*)/C_A\cong(\Z,+)$.
Let $C_B=H^*\cap (L_2 \times...\times L_{n+1})$. As $H^*\cap L_1=C_A$, we see that  $C_B$ is a disjoint union of $\{\gamma_1 H_B,...,\gamma_k H_B\}$ and
\[
H^*=\bigcup_{m\in \Z}(\eta^m C_A)\times (f(\eta^m)C_B)
\]

Since $\pi_{(2,...,n+1)}(H^*)/C_B\cong \Z$, applying $\MII_n$, there is a definable map  $f^*:\pi_1(H^*)\rightarrow \pi_{(2,...,n+1)}(H^*)$  such that
\[
 \bar f^*: \pi_1(H^*)/C_A\rightarrow \pi_{(2,...,n+1)}(H^*)/C_B,\ \ a/C_A\mapsto f^*(a)/C_B
\]
is an isomorphism, so $f^*(x)/C_B=f(x)/C_B$ as $(\Z,+)$ has only one automorphism. So
\[
H^*=\bigcup_{a\in \pi_1(H)}(a C_A)\times (f^*(a)C_B)
\]
By $\MI_n$, $C_A$ and $C_B$ are definable, so $H^*$ is definable.

Let $S_A=\bigcup_{m\in \Z}\eta^m H_A$ and $S_B=\bigcup_{m\in \Z}f(\eta)^m H_B$. Then $S_A$ and $S_B$ are definable by $\MI_n$. For any $(x ,y)\in H^*$, we have
\[
(x,y)\in H \iff (\bigwedge_{i=1}^k ({\theta_i}^{-1} x\in S_A\leftrightarrow {\gamma_i}^{-1} y\in S_B)),
 \]
so $H$ is definable as required.
\end{proof}

It is easy to see from  Lemma \ref{open subgp of Gm}, Lemma \ref{MII1}, Lemma \ref{MIIn+1}, and Lemma \ref{MIn+1} that Theorem \ref{Multi-Open Subgroup Theorem I} holds, namely, every open subgroup of $\Q_{p}^{\times}\times...\times \Q_{p}^{\times}$ is definable.

\begin{Lemma}\label{Quotient-group-of-multi-group}
Suppose that $G=L_1\times ...\times L_n$ with each $L_i=\Q_{p}^{\times}$. Let  $O\leq H$ be open subgroups of $G$.    Then $H/O$ is either finite or contains a torsion-free element (i.e. an element of infinite order).
\end{Lemma}
\begin{proof}
Induction on $n$.

Suppose that $H/O$ is infinite. Let $H_A=H\cap L_1$, $O_A=O\cap L_1$, $H_B=H\cap (L_2\times ...\times L_n)$, and $O_B=O\cap (L_2\times ...\times L_n)$.

\vspace{0.4cm}
\noindent $\bullet$ Case 1. $\pi_1(H)/\pi_1(O)$ is infinite. Then $\pi_1(O)\subseteq \Zp^\times$ and $\pi_1(H)$ is non-compact. Take any $a\in \pi_1(H)\backslash \Zp^\times$, then $(a,f_H(a))/O $ is torsion-free in $H/O$.

\vspace{0.4cm}
\noindent $\bullet$ Case 2. $\pi_1(H)/\pi_1(O)$ is  finite and   $H_A/O_A$ is infinite.   Then $O_A\subseteq \Zp^\times$ and $H_A$ is non-compact. Take any $a\in H_A\backslash \Zp^\times$, then $(a,1,...,1)/O $ is torsion-free in $H/O$.

\vspace{0.4cm}
\noindent $\bullet$ Case 3. $\pi_1(H)/\pi_1(O)$   and $H_A/O_A$ are finite.  Let $H^*=H\cap (\pi_1(O)\times L_2\times...\times L_n)$. Then $H/H^*\cong\pi_1(H)/\pi_1(O) $ is finite. So $H^*/O$ is infinite. To see $H/O$ has a torsion-free element, it suffices to find a torsion-free element in $H^*/O$. Without loss of generality, we assume that $\pi_1(H)=\pi_1(O)$.

Now $H_AO/O\cong H_A/O_A$ is a finite subgroup of $H/O$. Let $H_AO=O^*$, then $O^*\cap L_1=H_A$ and $H/O^*$ is an infinite quotient group of $H/O$. To see $H/O$ has a torsion-free element, it suffices to find a torsion-free element in $H/O^*$. Without loss of generality, we assume that $O_A=H_A$.

So we only need to consider the case where $\pi_1(H)=\pi_1(O)$ and $O_A=H_A$. Which means that $H/O\cong H_B/O_B$, and hence has a torsion-free element by our induction hypothesis.

\end{proof}

\begin{Lemma}\label{Lemma-open-subgp-II}
Suppose that $G = A\times B$ where $A$ is a (finite) product of copies of $(\Q_{p},+)$ and   $B$ is a (finite)  product of copies of $\Q_{p}^{\times}$.  If $H$ is an open subgroup of $G$. Then $H$ is definable.
\end{Lemma}
\begin{proof}
Let $\pi_A$ and $\pi_B$ be the projections from $G$ onto $A$ and $B$ respectively. Then  $H_A=G\cap A$  and  $H_B=G\cap B$  are definable by Theorem \ref{Additive-Open-Subgroup-is-definable} and Theorem \ref{Multi-Open Subgroup Theorem I}. By Lemma \ref{Quotient-group-of-add-group} and Lemma \ref{Quotient-group-of-multi-group}, $\pi_A(H)/H_A$ is a torsion group and $\pi_B(H)/H_B$ is either finite or has a torsion free element.  We conclude that  $\pi_A(H)/H_A\cong \pi_B(H)/H_B$ is finite. So $H$ is a finite union of translates of  $H_A\times H_B$, thus is definable.

\end{proof}

Note that Lemma \ref{Lemma-open-subgp-II} is precisely Proposition 2.2, which as noted in Section 2 suffices to prove the main Theorem 1.3 that every open subgroup of the group of $\Q_{p}$-points of a connected commutative algebraic group over $\Q_{p}$ is $p$-adic semialgebraic.



\section{The reductive case}
By a reductive algebraic group over $\Q_{p}$ we mean a connected algebraic group over $\Q_{p}$ which is an almost direct product of a semisimple group $G$ and an algebraic torus $T$, both of which will be over $\Q_{p}$.

We will answer Question 1.1 positively when $G$ is a reductive algebraic group.   In fact we will obtain a more general statement. Recall first that if $G$ is a connected algebraic group over $k$ then $G$ has a unique normal connected linear subgroup $L$, $L$ is defined over $k$, and $G/L$ is an abelian variety over $k$.

\begin{Proposition} Let $G$ be a connected algebraic group over $\Q_{p}$, such that the linear part $L$ of $G$ is an almost direct product of a semisimple group (over $\Q_p$) and a connected commutative linear algebraic group (over $\Q_{p}$). Then any open subgroup of $G(\Q_{p})$ is $p$-adic semialgebraic.
\end{Proposition}

We work towards the proof of Proposition 5.1. We will be referring to  \cite{Borel-Tits} and \cite{Platonov-Rapinchuk}.   We will let $k$ denote $\Q_{p}$. Fix a semisimple algebraic group $G$ over $k$, which we recall means that $G$ is connected and its  solvable radical (largest connected normal subgroup of $G$, necessarily defined over $k$) is trivial.  $G$ is said to be almost $k$-simple if it has no nontrivial normal connected subgroup defined over $k$. Sometimes one just says $k$-simple in place of almost $k$-simple, and this is what we will do.

\begin{Fact}  Assume $G$ is semisimple and over $k$ then $G$ is an almost direct product of  $k$-simple groups defined over $k$.
\end{Fact}

Let $G$ be semisimple and defined over $k$. $G$ is said to be $k$-isotropic.
 if $G$ has a nontrivial subgroup $T$ defined over $k$ which is a $k$-split torus, namely $T$ is isomorphic over $k$ to a product of copies of the multiplicative group. Otherwise $G$ is said to be $k$-anisotropic. See Theorem 3.1 of \cite{Platonov-Rapinchuk} for the following:

\begin{Fact} If $G$ is semisimple, over $k$ and $k$-anisotropic, then $G(k)$ is compact.
\end{Fact}

Recall that a connected algebraic group $G$ is said to be {\em simply connected} if whenever $H$ is a connected algebraic group and $f:H\to G$ is a surjective homomorphism of algebraic groups with finite (central) kernel then $f$ is an isomorphism.

\begin{Fact} Assume $G$ to be semisimple and over $k$. Then there is a simply connected group $\tilde G$ which is also over $k$ and a surjective homomorphism $f:{\tilde G}\to G$ over $k$ with finite (central) kernel.
\end{Fact}

See Proposition 2.10 of \cite{Platonov-Rapinchuk}.

\begin{Definition} Let $G$ be a semisimple group over $k$. By $G(k)^{+}$ we mean the (abstract) subgroup of $G(k)$ generated by its unipotent elements.
\end{Definition}

The following is due to Tits \cite{Tits}:
\begin{Proposition} Let $G$ be $k$-simple. Then  every subgroup of $G(k)$ normalized by $G(k)^{+}$ is either central in $G$ or contains $G(k)^{+}$. In particular $G(k)^{+}$ is abstractly simple, modulo its finite centre.

\end{Proposition}

The following is the truth of  the Kneser-Tits conjecture for $k = \Q_{p}$.  See Theorem 7.6 of \cite{Platonov-Rapinchuk}.
\begin{Proposition} Suppose $G$ is $k$-simple, isotropic, and simply connected. Then $G(k)^{+} = G(k)$.
\end{Proposition}

Here is a version (or consequence)  of the  unpublished theorem of Tits the proof of which, by Prasad, is in  \cite{Prasad}.

\begin{Theorem} Let $G$ be semisimple, $k$-simple (over $k$) and $k$-isotropic. Then
\newline
(i) $G(k)^{+}$ is semialgebraic and of finite index.
\newline
(ii) Any open subgroup $H$ of $G(k)$ is either compact (so semialgebraic), or contains $G(k)^{+}$, so of finite index and semialgebraic.
\end{Theorem}
\begin{proof}  Let $f:{ \tilde G} \to G$ be the ``universal cover" of $G$ given by Fact 5.4. By Proposition 5.7,  ${\tilde G}(k) = {\tilde G}(k)^{+}$.
By Corollary 3.20 of \cite{Borel-Tits}, $f({\tilde G}(k))$ (a semialgebraic subgroup of $G(k)$) is an open normal subgroup of $G(k)$ of finite index.  But $G(k)^{+} = f({\tilde G}(k)^{+}) = f({\tilde G}(k))$, so we have established that $G(k)^{+}$ is a semialgebraic (open, normal) subgroup of $G(k)$ of finite index.
Now let $H$ be an open subgroup of $G(k)$.  So $H\cap G(k)^{+}$ is an open subgroup of $G(k)^{+}$.  The result ``Theorem (T)"  of \cite{Prasad}  says that any proper open subgroup of $G(k)^{+}$ is compact. So either $H\supseteq G(k)^{+}$ in which case $H$ is semialgebraic of finite index in $G(k)$, or $H\cap G(k)^{+}$ is compact and open in $G(k)$, so semialgebraic, and as $G(k)^{+}$ has finiite index in $G(k)$, $H \cap G(k)^{+}$ has finite index in $H$ so the latter is also compact and semialgebraic.
\end{proof}

\begin{Lemma}  Let $G$ be a connected algebraic group over $k$ which is an almost direct product of connected algebraic groups $G_{1}$ and $G_{2}$ over $k$ where $G_{1}$ is semisimple and $k$-simple. Let $H$ be an open subgroup of $G(k)$, let $A = H\cap G_{1}(k)$ and $B = H\cap G_{2}(k)$. Then the group $AB$ has finite index in $H$.
\end{Lemma}
\begin{proof} There is no harm in assuming that $G$ is a direct product $G_{1}\times G_{2}$  (with $G_{1}$ being $k$-simple).
Let $A_{1} = \pi_{1}(H)$, the projection of $H$ on $G_{1}$.  Then $A\subseteq A_{1}$ are open subgroups of $G_{1}$ with $A$ normal  in $A_{1}$.

Note that $H/(A\times B)$ is naturally isomorphic to $A_{1}/A$, so it suffices to show that $A$ has finite index in $A_{1}$.
Now if $G_{1}$ is $k$-anisotropic then $G_{1}(k)$ is compact, so both $A$ and $A_{1}$ have finite index in $G_{1}(k)$, in particular $A_{1}/A$ is finite.
Otherwise $G_{1}$ is   $k$-isotropic.  Consider $A$. We use Theorem 5.8.  If $A$ is compact, then $A\cap G(k)^{+}$ has finite index in $A$.  Now if $A_{1}$ contains $G(k)^{+}$ then $A\cap G(k)^{+}$ is a proper normal (infinite) subgroup of $G(k)^{+}$ contradicting Proposition 5.6.  So $A_{1}$ is compact.  This $A_{1}/A$ is finite.
Otherwise $A$ has  finite index in  $G(k)$ so also in $A_{1}$.
\end{proof}

\vspace{5mm}
\noindent
{\em Proof of Proposition 5.1.} As in the proof of Lemma 2.1, we may assume that $G$ is linear (as if $A$ is the maximal abelian variety quotient of $G$ then  $A(k)$ is compact).
Let $G$ be the almost direct product of a semisimple group $S$ (over $k$) and a commutative group $D$ over $k$ (both connected).
By Fact 5.2 we can write $S$ as the almost direct product of $k$-simple groups $G_{1},..,G_{r}$. We do induction on $r$. If $r=0$ then $G = D$ is commutative and use Theorem 1.3.
Suppose $r = 1$, so $G$ is the almost direct product of $G_{1}$ and $D$. Let $H$ be an open subgroup of $G$. Let $A = H\cap G_{1}(k)$ and $B =H \cap D(k)$. Then $A$ is semialgebraic  if $G_{1}$ is $k$-isotropic (Theorem 5.8), and compact, so semialgebraic when $G_{1}$ is $k$-anistropic.  $B$ is semialgebraic by Theorem 1.3. So $AB$ is semialgebraic and by Lemma 5.9 so is $H$.

Now suppose $r > 1$. Write $G$ as the almost direct product of $G_{1}$ and a group $L$ where $L$ is an  almost direct product of $G_{2},...,G_{r}, D$.
Let $H$ be an open subgroup of $G$, and let $A= H\cap G_{1}(k)$, and $B = H\cap L(k)$.  By induction, both $A$ and $B$ are semialgebraic, hence so is $AB$ and so also $H$, by Lemma 5.9.
\qed


\begin{thebibliography}{99}

\bibitem{Borel-Tits} A. Borel and J. Tits, Homomorphismes "abstraits" de groupes algebriques simples, Annals of  Mathematics,  vol 97 (1973), 499 - 571.


\bibitem{Macintyre} A. Macintyre, On definable subsets of $p$-adic fields, JSL, 41 (1976), 605 - 610.

\bibitem{Milne} J. S. Milne, Algebraic groups, the theory of group schemes of finite type over a field, Cambridge Studies in advanced mathematics, vol. 170, CUP, 2017.

\bibitem{NP}  A. Nesin and A. Pillay,  Open subgroups of $GL_{2}(\Q_{p})$, Seminarbericht nr. 104, Proceedings of 7th Easter Conference on Model Theory, (editors Dahn, Wolter), Sektion Math. of Humboldt Univ. Berlin, Deutsche Demokratische Republik, Berlin, 1989.

\bibitem{Onshuus-Pillay}  A. Onshuus and A. Pillay, Definable groups and compact $p$-adic Lie groups, J. London Math. Soc. 78 (2008), 233 -247.


\bibitem{Pillay}  A. Pillay, An application of model theory to real and $p$-adic algebraic groups, J. Algebra 126 (1989), 139-146.


\bibitem{Platonov-Rapinchuk} V. Platonov and A. Rapinchuk, Algebraic Groups and Number Theory,  Academic Press, 1994.


\bibitem{Prasad} G. Prasad, Elementary proof of a theorem of Bruhat-Tits-Rousseau and of a theorem of Tits, Bull Soc. math. France, 110 (1982), 197 - 202.




\bibitem{definable Skolem functions} P. Scowcroft, A note on definable Skolem functions, J. Symbolic Logic, 53 (1988), 905-911.




\bibitem{Tits} J. Tits, Algebraic and abstract simple groups,  Annals of Mathematicas, vol. 80 (1964), 313 - 329.







\end{thebibliography}
\end{document}